\documentclass[twoside,leqno,11pt]{article}
\usepackage{alltt}
\usepackage{psfig}

\newcommand{\email}[1]{{\sf #1}}
\newenvironment{CE}{%
\begin{small}\begin{alltt}}{%
\end{alltt}\end{small}}

\newcommand{\C}[1]{\begin{tt}#1\end{tt}}

\newcommand{\refsec}[1]{Section~\ref{#1}}

\begin{document}
\cleardoublepage
\pagestyle{headings}

\title{Active Libraries: Rethinking the roles of compilers and libraries}

\author{Todd L. Veldhuizen and Dennis Gannon\thanks{Indiana University Computer Science 
Department.  
\email{tveldhui@extreme.indiana.edu}}}

\date{}
\maketitle

\pagenumbering{arabic}

\begin{abstract}
\noindent
We describe {\it Active Libraries}, which take an active role
in compilation.  Unlike traditional libraries which are passive
collections of functions and objects, Active Libraries may generate components,
specialize algorithms, optimize code, configure and tune themselves
for a target machine, and describe themselves to tools (such as
profilers and debuggers) in an intelligible way.  Several 
such libraries are described, as are implementation technologies.
\end{abstract}

\section{Introduction}

This paper attempts to document a trend 
toward libraries which take an active role in generating
code and interacting with programming tools.
We call these {\em Active Libraries}.
They solve the problem of how to provide
efficient domain-specific abstractions
(\refsec{s:whyal}): 
active libraries are able to define abstractions,
and also control how they are optimized.
Several existing libraries for
arrays, parallel physics, linear algebra and
Fast Fourier Transforms fit the description of
active libraries (\refsec{s:examples}).
These libraries take novel approaches to
generating optimized code (\refsec{s:optimization}).
To implement active libraries, programming
systems and tools which open up the development
environment are needed (\refsec{s:technologies}).

\subsection{Why are active libraries needed?}

\label{s:whyal}

\noindent
To produce readable, maintainable scientific computing codes,
we need abstractions.  
Every subdomain in scientific computing has its
own requirements: interval arithmetic, tensors, polynomials,
automatic differentiation, sparse arrays, spinors, meshes, and so on.  
How should these abstractions be provided?  The alternatives are:

{\bf Extend mainstream languages.}  
In the past, syntax and efficiency concerns have encouraged
building abstractions into languages and compilers.  Fortran 95, 
for example, has built-in complex numbers and numerical arrays,
and many intrinsic functions for common numerical operations.
Efforts are underway to extend the Java language for
scientific computing \cite{JavaGrande}, and similar efforts
in the past have extended C \cite{Plauger:1993:NCE}.
However, loading up these languages with features
for the scientific market can meet with only limited
success.
The costs of compiler development
are quite high, and scientific computing is
a comparatively small segment of the market.
Even Fortran, whose bread and butter is scientific
computing, is showing signs of hitting economic limits
to its size.  A controversial 60-page proposal to add
interval arithmetic to Fortran 2000 was eventually discarded
after much debate.  The committee had to balance the
limited demand for interval arithmetic against the large
implementation costs for vendors.  Although we may
succeed in getting mainstream languages to incorporate basic
scientific features such as numeric arrays and complex numbers,
the likelihood of these languages providing 
specialized features such as interval arithmetic and sparse 
arrays is small.

There are other disadvantages to building abstractions
into mainstream languages:  feature turnaround is very slow,
since extensions require championing a proposal through years
of standards committee meetings.  Experimental features must
be prematurely standardized, before experience can produce
a consensus on the ``right way'' to implement them.


{\bf Domain-Specific Languages (DSLs).}
Dozens of DSLs for scientific computing have been produced
to handle sparse arrays,
automatic differentiation, interval arithmetic, adaptive mesh
refinement, and so on.  DSLs are often implemented
as preprocessors for mainstream languages such as
C, Fortran, C++, or Java.
The growing availability of compiler
construction tools has encouraged a proliferation of DSLs.
Although DSLs are an attractive alternative,
many people would prefer to work in mainstream languages 
for the wealth of tools support and
libraries available.  DSLs frequently have problems with
portability and long-term support, since they tend to
be research projects.  It is also impossible to use
multiple DSLs at once; for example, although separate 
Fortran-based DSLs are
available for both sparse arrays and interval arithmetic,
you cannot use the features of both DSLs in the same
source file.

{\bf Object-oriented language features.}
An alternative to building abstractions into languages
is to provide language features which allow library developers
to construct their own abstractions.  In C++ and Fortran 90/95, we
are seeing libraries for many applications (sparse arrays,
interval arithmetic, data-parallel arrays)
which were previously solved by domain-specific languages.
Unfortunately, such libraries are hard to optimize.  Compilers
have difficulty because they lack semantic knowledge of the
abstractions: instead of seeing array operations, they see
loops and pointers.  Libraries also tend to have layers of
abstraction and side effects which confound optimizations.
The heroic optimizers needed to overcome these problems may
never appear, because the economics of the scientific market
may not support their development.  
It is doubtful that the optimization problems admit
a general-purpose solution, since every problem domain
has its own tricks and peculiarities.

{\bf What we really need} are language features
which allow library developers to define their own
abstractions, and also to specify how these abstractions
are optimized.
We call this solution {\it Active Libraries.}
Active Libraries combine the benefits of built-in 
language abstractions (nice syntax and efficient code) with those 
of library-level abstractions
(adaptability, quick feature turnaround, cheap to implement).

\section{Examples of Active Libraries}

\label{s:examples}

In defining Active Libraries, we are not proposing a new
concept, but rather trying to summarize what many people
are already trying to do.  In the following sections,
we highlight existing software packages which illustrate the
characteristics of Active Libraries.

\subsection{Blitz++}

\label{s:Blitz}

The Blitz++ library \cite{Veldhuizen98a} provides
generic array objects for C++ similar to those in Fortran 90,
but with many additional features.
In the past, C++ array libraries have been 3-10 times
slower than Fortran, due to the temporary arrays
which result from overloaded operators.
%
%
Blitz++ solves this problem using the {\it expression templates}
technique \cite{Veldhuizen95b} to
generate custom evaluation kernels for array expressions.
The library performs many loop transformations (tiling,
reordering, collapsing, unit stride optimizations, etc.) which
have until now been the responsibility of optimizing compilers.
Blitz++ also generates different code depending on the
target architecture.  

For operations on small vectors and matrices, Blitz++
uses the {\it template metaprogram} technique
\cite{Veldhuizen95a} to generate specialized algorithms.
This avoids the performance penalty often associated
with small objects by completely unrolling loops and
inlining code.

\subsection{POOMA}

\label{s:POOMA}

POOMA is a C++ library for parallel physics which uses many of
the same techniques as Blitz++.
Users write simple array expressions, such as ``A=B+C'',
which trigger the generation of data-parallel implementation
routines using threads and message passing \cite{POOMA}.
POOMA uses template techniques to generate components (such as
Fields) using a variety of types, geometries, addressing
schemes, data distribution and communication parameters.

\subsection{Matrix Template Library}

\label{s:MTL}
The Matrix Template Library (MTL) \cite{blais_poosc} is a C++
library which extends the ideas of STL
\cite{Musser94} to linear algebra.  MTL handles both sparse
and dense matrices.
For dense matrices, MTL uses template metaprograms to
generate tiled algorithms.  Tiling is a crucial technique
for obtaining top performance from cache-based memory
systems; MTL uses template metaprograms to tile on both
the register and cache level.  For register tiling,
it uses template metaprograms to completely unroll loops.
MTL provides
generic, high-performance algorithms which are competitive
with vendor-supplied kernels.

\subsection{Generative Matrix Computation Library}

The Generative Matrix Computation Library (GMCL)
\cite{GMCL} provides heavily parameterized
matrix classes.  Users can specify the element
type, whether the matrix is dense or sparse, the
storage format (including several sparse formats),
dynamic or static memory allocation, error checking,
and several other parameters.
The GMCL uses template metaprograms to examine the
parameters and determine any interactions between
them (for example, sparse matrices cannot use
static memory allocation); it then instantiates
a matrix class with the desired characteristics.
The implementation is roughly 7500 lines of C++
code, yet covers more than 1840 different kinds
of matrices.  Despite this flexibility, the authors
report performance on par with manually generated
code.


\subsection{FFTW: The Fastest Fourier Transform in the West}

\label{s:FFTW}

FFTW \cite{Frigo97} is a C library for Fast Fourier Transforms.
It generates a collection of small 
``codelets,'' each of which is a small step in the transform
process.  At installation, FFTW evaluates the performance
of each codelet for the target architecture.  At run-time,
the codelets are dynamically stitched together to
perform FFTs.  FFTW records the ordering of these
codelets using bytecode, which is generated and
interpreted at run time.
Based on extensive benchmarking, the authors
report superior performance over all commonly used FFT
packages.

\subsection{PhiPAC, ATLAS}

\label{s:PhiPAC}

Obtaining top performance for matrix operations requires
substantial expertise and hand-tuning.  PhiPAC \cite{bilmes97b}
provides a methodology to achieve near-peak performance automatically.
It uses parameterized code generators, whose parameters are
related to machine-specific tuning.  PhiPAC searches
the parameter space to find the best implementation
for a given target architecture.  On dense matrix-matrix
multiplication, PhiPAC performs better then vendor-supplied
kernels on many platforms.  The ATLAS package \cite{Whaley98} provides
similar capabilities.

\section{Optimization models}

\label{s:optimization}

Active Libraries make it possible to approach
code optimization in several new and exciting ways.
In describing some of the new approaches to optimization
being explored, it is useful to distinguish 
between two flavours of optimization, which we call
{\it low-level} and {\it high-level}:

\begin{itemize}
\item {\bf Low-level} optimizations can be applied without
knowing what the code is supposed to do (copy
propagation, dead code elimination, instruction
scheduling, loop pipelining).

\item {\bf High-level} optimizations which require some understanding
of the operation being performed.  Examples: tiling
for stencils; fusing loops over sparse arrays;
iteration-space tiling.
\end{itemize}

\subsection{Transformational optimization}

Traditional approaches to optimization are {\it transformational}.
In transformational optimization, low-level code is transformed
into an equivalent (but hopefully faster) implementation.
Typical transformational optimizations are constant propagation,
dead code elimination, and loop transformations.

One of the difficulties with transformational optimization is
that the optimizer lacks an understanding of the {\it intent}
of the code.  This makes it harder to apply high-level
optimizations.  For example, rather than seeing an array stencil
operation, a transformational optimizer will
just see loops, pointers, and variables.  To apply
interesting optimizations, the optimizer must
recognize that this low-level code represents
a stencil operation (or more accurately, that it
possesses the sort of data dependencies which benefit
from tiling).  More generally, the optimizer
must infer the {\it intent} of the code to apply
higher-level optimizations.  To this end,
sophisticated optimizers employ algebraic
reasoning, pattern
recognition and matching techniques.
However, such optimizers can still only apply
radical optimizations in simple situations, and
with limited success.

Another problem with transformational optimizers is
the lack of extensibility.  If the optimizer doesn't
recognize what your code is trying to do, you are
probably out of luck.  Optimizations for dense
arrays are fairly reliable, because their use is
very common; if you are working
with sparse arrays or interval arithmetic, the
likelihood of achieving optimal performance is
smaller.


\subsection{Generative optimization}

Many environments which provide higher-level abstractions
(for example, arrays in Fortran 90) can use {\it generative
optimization}.  Code written using higher-level abstractions
can be regarded as a
specification of what operation needs to be performed;
an efficient implementation is then generated to
fulfill the ``specification''.  This approach to
optimization can be much simpler than transformational
optimization, since the full semantics are available
to whatever generates the optimized code.  Generative
optimization can make it simpler to apply radical,
high-level optimizations.

While generative optimization can produce good
code for individual operations, it tends to miss
optimizations which depend on context.
For example, one problem encountered in libraries
which use expression templates (e.g. Blitz++, POOMA)
is that while individual array statements can
be optimized well, opportunities for between-statement
optimizations are missed.  For example, in the
array statements \C{A=B+C; D=B-C;} there is
a substantial gain if the two expressions
are evaluated simultaneously, since they share
the same operands.  Attempts to solve this
problem have focused on increasing granularity,
so that code is generated
for basic blocks of array statements, rather
than individual statements.

\subsection{Explorative optimization}

For library developers, finding a near-optimal
implementation of a routine is very difficult.
Modern architectures can behave unpredictably due
to pipeline and cache effects.  Between the
library writer and the hardware lies the
compiler, a black box which transforms one's
code in sometimes mysterious ways.  Aside from some
basic guidelines relating to cache reuse,
performance tuning usually requires
randomly adjusting code and measuring the
result.

{\it Explorative optimization} gives up on
the notion that performance is somehow
predictable.  The basic approach is to
examine an algorithm and identify parameters
which might affect performance (for
example loop structures, tile sizes, 
and unrolling factors).  One then writes
a parameterized code generator which
produces variants of the basic algorithm.
The parameter space is then explored point
by point, and the performance of each
variant is measured until the best implementation
is found.  This approach was pioneered
by FFTW (\refsec{s:FFTW}) and
PhiPAC (\refsec{s:PhiPAC}).
On dense matrix-matrix
multiplication, PhiPAC performs better
than vendor-supplied kernels on many
platforms.  

Explorative optimization is expensive in time,
but is worth it if a near peak-performance
kernel is required.  So far, there is
no way to automatically construct variants
of a given piece of code; one must
write the code generators manually.


\subsection{Compositional optimization}

Compositional optimization is useful when
a problem can be decomposed into a sequence
of calls to well-tuned kernels.  The
FFTW package (\refsec{s:FFTW}) uses this approach
to decompose FFTs into a sequence of calls
to high-speed kernels which were found using
explorative optimization.  FFTW
records and interprets the sequence of calls
using a bytecode.
Compositional optimization has also been
used in the context of compilers, to generate
high-performance communication code
\cite{Stichnoth97}.

For compositional optimization to be
effective, the problem must be decomposable
into coarse chunks, so that the overhead
of composition is negligible.

\section{Technologies for Active Libraries}

\label{s:technologies}

In the following sections, we describe technologies relevant
to Active Libraries.  Some of these are already being
used, for example generic programming and C++ templates;
others are promising technologies still being explored.

\subsection{Component generation}

Traditional libraries define {\it concrete components}.
We use the term {\em components} loosely here, to mean
an algorithm, class, or collection of these things.
By concrete, we mean that the behavior is fixed: the
component operates on a fixed kind of data (for example,
arrays of \C{double}), using a specific data
structure, and handles errors in a fixed way.
Concrete components might allow for some flexibility
through configuration variables,
but the code of the component is unchanging.
Concrete components require a trade-off between
flexibility and efficiency: 
if the library developer
wants to provide a customizable component, this
must be done through callbacks and runtime checking of
configuration variables, which are often inefficient.
To solve this problem, we need ways to {\em generate}
customized components on demand.
We contrast two approaches: Generic Programming 
and Generative Programming.


\subsubsection{Generic Programming}

The aim of Generic Programming can be summarized
as ``reuse through parameterization''.  Generic components
have parameters which customize their behaviour.  When a 
generic component is {\em instantiated} using a particular
choice of parameters, a concrete component 
is generated.
This allows library developers to create components which
are very customizable, yet retain
the efficiency of statically configured code.
Probably the greatest achievement of the Standard
Template Library in C++ \cite{Musser94} was to
separate algorithms from the data structures on which
they operate, allowing them to be combined in a mostly
orthogonal way.  In C++, template parameters can 
effectively be types, data structures, or even pieces 
of code and algorithms.  Other languages which support
generic programming are Ada (via its {\it generics}
mechanism), and Fortran 2000 (albeit in a limited way).

There are two main benefits of generic programming:
(1) Library developers have to develop and maintain less
code; (2) Application developers find it easier to
find components which match their needs.
There are some limitations of generic programming
as currently supported:  (1) Aggressive parameterization
is quite ugly; providing more than 2-3 template parameters
introduces serious usability problems.  This can be
alleviated somewhat by named template parameters.\footnote{Named 
template parameters can be listed in any order, with missing
parameters assuming default values.  Thought not directly supported
in C++, the technique can be faked (see \cite{GMCL}).}
(2) Generic programming may result in code bloat, due to
multiple versions of generic components.
(3) Perhaps most importantly, generic programming
limits code generation to substituting concrete
types for generic type parameters, and welding together 
pre-existing fragments of code in fixed patterns.  It does not 
allow generation of completely new code, nor does it
allow computations
to be performed at compile time.\footnote{We regard techniques such as
template metaprograms and expression templates to lie
outside the domain of generic programming.}

\subsubsection{Generative Programming}

Generative Programming \cite{Czarnecki98b,KrUl99}
is a broader term which encompasses generic programming, 
code generation, code analysis and transformation, 
and compile-time computation.  In general, it refers
to systems which generate customized components to
fulfill specified requirements.  

A central idea in Generative Programming (and also
in Aspect-Oriented Programming \cite{AOP}) is {\it separation 
of concerns}: the notion that important issues should be
dealt with one at a time.
In current languages, there are many aspects
such as error handling, data distribution, and
synchronization which cannot be dealt with in a localized way.  
Instead, these aspects are scattered throughout the code.
One of the goals of Generative Programming is to 
separate these aspects into distinct pieces of code.
These pieces of code are combined to produce a
needed component.  Doing so often requires more than
cutting and pasting, and this is where the need for
code generation, analysis, and transformation arises.
Active Libraries can be viewed as a vehicle for implementing the
goals of Generative Programming \cite{Czarnecki98b}.






\subsection{C++ Templates}

In addition to enabling generic programming, the template
mechanism of C++ unwittingly introduced powerful code
generation mechanisms.  Nested templates allow 
data structures to be created and manipulated at
compile time, by encoding them as types. This is
the basis of the {\it expression templates} technique~\cite{Veldhuizen95b},
which creates parse trees of array expressions at compile time,
and uses these parse trees to generate customized kernels.
Expression templates also provides a crude facility similar
to the {\tt lambda} operator of functional languages, which
may be used to replace callbacks.  
{\em Template metaprograms}~\cite{Veldhuizen95a} 
use the template instantiation mechanism to perform 
computations at compile-time, and generate specialized algorithms
by selectively inlining code as they are executed.
Although these techniques are powerful, the accidental
nature of their presence has resulted in a clumsy syntax.
Nonetheless, they provide a reasonable way to implement
Active Libraries in C++, and several libraries based on
these techniques (Blitz++, POOMA, MTL) are being
distributed.  Recently, a package which simplifies
the construction of template metaprograms has been
made available \cite{CE98}.


\subsection{Extensible compilation, Reflection, and Metalevel Processing}

In metalevel processing systems, library writers are
given the ability to directly manipulate language constructs.
They can analyze and transform syntax trees, and generate
new source code at compile time.  The MPC++ metalevel
architecture system \cite{Ishikawa96}
provides this capability for the C++ language.  MPC++
even allows library developers to extend the syntax
of the language in certain ways (for example, adding
new keywords).  Other examples of metalevel processing
systems are Xroma \cite{Czarnecki98b}, Open C++ \cite{OpenC++}, 
and Magik \cite{EnglerDSL}.
A potential disadvantage of metalevel processing systems is the
complexity of code which one must write: modern languages have
complicated syntax trees, and so code which manipulates
these trees tends to be complex as well.



\subsection{Run-Time Code Generation (RTCG)}


RTCG systems allow libraries to generate customized code
at run-time.  This makes it possible to perform optimizations
which depend on information not available until run-time,
for example, the structure of a sparse matrix or the
number of processors in a parallel application.
Examples of such systems which generate native
code are `C (Tick-C) \cite{POPL96*131,EnglerPLDI},
and Fabius \cite{LeoneMarka1994a}.  Speeds as high as
6 cycles per generated instruction have been achieved.
Recently, this technology has been extended to
C++ \cite{Fujinami97}.


\subsection{Partial Evaluation}

Code generation is an essential part of active libraries.  
Over the past two decades, researchers in the field of
Partial Evaluation have developed an extensive theory
and literature of code generation \cite{Jones:1996:IPE}.  
In its simplest form, a partial evaluator analyzes a program and
separates its data into a static portion (values known
at compile-time) and a dynamic portion (values not
known until run-time).  It then evaluates as much
of the program as possible (using the static values)
and outputs a specialized {\em residual} program.
For example, a partial evaluator could take a
dot-product routine, and produce a specialized
version for a particular vector length.  These
techniques have been applied to scientific codes
with promising results \cite{AIM-1487,Berlin90,KKZG:95:FortranSpec}.

However, this just provides a taste of the field; 
partial evaluation has evolved into a comprehensive toolbox containing
both theories and practical software.  One of the most important
theoretical contributions was that the concept of
{\em generating extensions}~\cite{Ershov:78} unifies
a very wide category of apparently different program
generators.  Using partial evaluation, concrete components
which check configuration variables at run-time
can be transformed into {\em component generators}
(or {\em generating extensions} in the terminology
of the field~\cite{BEJ:88,DGT:96,Jones:Gomard:Sestoft:93:PartialEvaluation})
which produce customized components (eliminating the
run-time checking of configuration variables).
Automatic tools for turning a general component
into a component generator now exist for C
and Scheme \cite{Jones:Gomard:Sestoft:93:PartialEvaluation}.
However, such tools are not yet available for C++.

\subsection{Multilevel Languages}

Another important contribution of Partial Evaluation
is the concept of two-level (or more generally, multi-level)
languages.  Two-level languages contain static constructs
(which are evaluated at compile-time) and dynamic code
(which is compiled and later evaluated at run-time).
Two-level languages provide a
simpler notation for writing code generators (compared
to systems which generate source code or intermediate
representations).  For example, consider 
a (fictional) two-level language based on C++, in which 
static variables and control flow are annotated
with the \C{@} symbol.  A code generator for dot products
of length N would be written as:

\begin{CE}
        double dot(double* a, double* b, int@ N) \{
            double sum = 0;
            for@ (int@ i=0; i < N; ++i)
                sum += a[i] * b[i];
            return sum;
        \}
\end{CE}

\noindent
In this example, the \C{@} symbol indicates that \C{N} and
\C{i} are compile-time variables.  The \C{for@} loop is
evaluated at compile time, and the residual code is equivalent
to \C{N} statements of the form \C{sum+=a[i]*b[i]}, with
\C{i} replaced by appropriate integer literals.
This dot-product generator is dramatically simpler than an equivalent
implementation using template metaprograms or metalevel
processing.
Two-level languages have been used to generate customized
run-time library code for parallel compilers
\cite{Stichnoth97}.








\subsection{Extensible Programming Tools}

Libraries typically have two (or more) layers: there
are user-level classes and functions, and behind them are
one or more implementation layers.  Using tools such as
debuggers and profilers with large libraries is troublesome,
since the tools make no distinction between user-level
and implementation code.  For example, a user who invokes
a debugger on an array class library will be confronted
with many irrelevant (and undocumented) private data members,
when all they wanted to see was the array data.  Using
a profiler with such a library will expose many implementation
routines, instead of indicating which array expressions were
responsible for a slow program.
This problem is compounded when
template libraries are used, with their long symbol names
and many instances.

The solution may be extensible tools, which provide
hooks for libraries to define customized support for
debugging, profiling, etc.  
An example of such interaction is Blitz++ and
the Tau \cite{Shende98} profiling package.
Tau is unique in that it allows libraries to
instrument themselves.  This allows Blitz++ to
hide implementation routines, and only expose
user-level routines.  Time spent in the library
internals is correctly attributed to the
responsible user-level routine.  
Blitz++ describes array evaluation kernels
to Tau using pretty-printing.  When users profile
applications, they do not see incomprehensible
expression template types; rather, they see expressions
such as ``A=B+C+D''.

If debuggers provided similar hooks, users could
debug scientific codes by looking at visualizations
and animations of the array data.  The general
goal of such tool/library interactions is to
provide a user-oriented view of libraries, rather
than an implementation-oriented view.


\section{Conclusions}


Active Libraries are able to define domain-specific
abstractions, and also control how these abstractions
are optimized.  This may involve compile-time
computations, code generation, and even code
analysis and transformation.  
It is no coincidence that all of the existing
examples are libraries for scientific
computing: this is a field which requires
many abstractions, and also demands
high performance.  Active libraries may
be the best way to construct
and deliver efficient,
configurable abstractions.

\section{Acknowledgements}

This work was supported in part by NSF grants CDA-9601632
and CCR-9527130.
Portions of this paper grew out of joint work with
David Vandevoorde (who coined the term {\it active library}),
Krzysztof Czarnecki, Ulrich Eisenecker, and Robert Gl{\"uck}
\cite{Czarnecki98b}.

\nocite{EnglerPLDI}

\bibliographystyle{siamproc}
\bibliography{activelib,robert,ulrich,partial,oon}

\end{document}